# Homological mirror symmetry for toric orbifolds of toric del Pezzo surfaces

Kazushi Ueda and Masahito Yamazaki

**Abstract**

We formulate a conjecture which describes the Fukaya category of an exact Lefschetz fibration defined by a Laurent polynomial in two variables in terms of a pair consisting of a consistent dimer model and a perfect matching on it. We prove this conjecture in some cases, and obtain homological mirror symmetry for quotient stacks of toric del Pezzo surfaces by finite subgroups of the torus as a corollary.

## 1 Introduction

The lattice of vanishing cycles equipped with the intersection form is called the *Milnor lattice*. It is a fundamental object in singularity theory, which is related to other fields of mathematics such as generalizations of root systems (see e.g. [Sai98]) and monodromy of hypergeometric functions. The Milnor lattice admits a categorification called the *Fukaya category of Lefschetz fibration*, defined by Seidel [Sei01b] based on an idea of Kontsevich [Kon98]. It is an $A_\infty$-category whose set of objects is a distinguished basis of vanishing cycles and whose spaces of morphisms are Lagrangian intersection Floer complexes. Although they are important invariants in singularity theory, it is often difficult to compute the Milnor lattice of a holomorphic function, let alone its Fukaya category.

Recent advances in string theory have given a significant progress in the case of Laurent polynomials in two variables. For a Laurent polynomial

$$W(x, y) = \sum_{(i,j) \in \mathbb{Z}^2} a_{ij} x^i y^j,$$

its Newton polygon $\Delta \subset \mathbb{R}^2$ is defined as the convex hull of $(i, j) \in \mathbb{Z}^2$ such that $a_{ij} \neq 0$;

$$\Delta := \mathrm{Conv}\{(i, j) \in \mathbb{Z}^2 \mid a_{ij} \neq 0\}.$$

We always assume that the origin is in the interior of $\Delta$. With a convex lattice polygon containing the origin in its interior, one can associate a pair $(G, D)$ consisting of a *consistent dimer model* $G$ and a *perfect matching* $D$ on $G$. Dimer models are introduced by string theorists to study supersymmetric gauge theory in four dimensions [FHM+06, FHV+06, FV06, HHV06, HK05, HV07]. The pair $(G, D)$ encodes the information of a directed $A_\infty$-category $\mathcal{A}^\rightarrow$ described in [FU10]. Motivated by Feng, He, Kennaway and Vafa [FHKV08], we formulate the following conjecture:

**Conjecture 1.1.** *For a convex lattice polygon $\Delta$ containing the origin in its interior, there exist*

- *a Laurent polynomial $W$ whose Newton polygon coincides with $\Delta$, and*
- *a pair $(G, D)$ consisting of a consistent dimer model and a perfect matching on it, whose characteristic polygon coincides with $\Delta$*

*such that one has a quasi-equivalence*

$$\mathcal{A}^{\rightarrow} \cong \mathfrak{Fuk}\, W \tag{1.1}$$

*of $A_\infty$-categories.*

With a convex lattice polygon $\Delta$ containing the origin in its interior, one can associate a two-dimensional toric Fano stack $X$. If $\Delta$ is the characteristic polygon of a pair $(G, D)$, then one has an equivalence

$$D^b \mathcal{A}^{\rightarrow} \cong D^b \operatorname{coh} X \tag{1.2}$$

of triangulated categories [IU, FU10]. By combining the equivalences (1.1) and (1.2), one obtains homological mirror symmetry

$$D^b \operatorname{coh} X \cong D^b \mathfrak{Fuk}\, W \tag{1.3}$$

conjectured by Kontsevich [Kon95, Kon98].

Now consider the pull-back

$$\widetilde{W} = W \circ \exp : \mathbb{C}^2 \to \mathbb{C}$$

of $W$ by the universal cover

$$\exp : \mathbb{C}^2 \to (\mathbb{C}^\times)^2$$

of the torus. The fact that $\widetilde{W}$ has infinitely many critical points for a given critical value does not cause any problem, and one can define the Fukaya category $\mathfrak{Fuk}\, \widetilde{W}$ of $\widetilde{W}$. We formulate Conjecture 6.3 in Section 6, which is slightly stronger than Conjecture 1.1 and implies a torus-equivariant version of homological mirror symmetry for two-dimensional toric Fano stacks:

**Theorem 1.2.** *If Conjecture 6.3 holds for a lattice polygon $\Delta$, then there is an equivalence*

$$D^b \operatorname{coh}^{\mathbb{T}} X \cong D^b \mathfrak{Fuk}\, \widetilde{W} \tag{1.4}$$

*of triangulated categories, where $\mathbb{T}$ is the algebraic torus acting on $X$ and $D^b \operatorname{coh}^{\mathbb{T}} X$ is the derived category of $\mathbb{T}$-equivariant coherent sheaves on $X$.*

Equivariant homological mirror symmetry (1.4) for $X$ implies homological mirror symmetry for the quotient stack $[X/A]$ of $X$ by a finite subgroup $A$ of the torus $\mathbb{T}$:

**Theorem 1.3.** *If Conjecture 6.3 holds for a lattice polygon $\Delta$, then homological mirror symmetry*

$$D^b \operatorname{coh} X \cong D^b \mathfrak{Fuk}\, W$$

*holds for any lattice polygon $\phi(\Delta)$ obtained from $\Delta$ by an integral linear transformation $\phi : \mathbb{Z}^2 \to \mathbb{Z}^2$ of rank two.*

We prove the following in this paper:

**Theorem 1.4.** *Conjecture 6.3 holds for lattice polygons corresponding to toric del Pezzo surfaces.*

This implies homological mirror symmetry for toric orbifolds of toric del Pezzo surfaces. The equivalence (1.3) is proved for $\mathbb{P}^2$ and $\mathbb{P}^1 \times \mathbb{P}^1$ by Seidel [Sei01a], weighted projective planes and Hirzebruch surfaces by Auroux, Katzarkov and Orlov [AKO08], and toric del Pezzo surfaces by Ueda [Ued06]. See also Abouzaid [Abo06, Abo09] for an application of tropical geometry to homological mirror symmetry, and Kerr [Ker08] for the behavior of homological mirror symmetry under weighted blowup of toric surfaces. Slightly different versions of homological mirror symmetry for toric stacks are proved by Fang, Liu, Treumann and Zaslow [Fan08, FLTZa, FLTZb] and by Futaki and Ueda [FU10]. Homological mirror symmetry for not necessarily toric del Pezzo surfaces is proved by Auroux, Katzarkov and Orlov [AKO06], and it is an interesting problem to extend Conjecture 6.1 to include non-toric cases. Since the $A_\infty$-category $\mathcal{A}^\rightarrow$ associated with a consistent dimer model and an internal perfect matching can be derived-equivalent only to toric stacks as in (1.2), one needs more general objects than dimer models for such an extension.

The organization of this paper is as follows: We collect basic definitions and facts on dimer models and $A_\infty$-categories in Section 2. Toric Fano stacks associated with lattice polygons are recalled in Section 3, and Fukaya categories of exact Lefschetz fibrations are recalled in Section 4. In Section 5, we discuss homological mirror symmetry for $\mathbb{P}^2$ proved by Seidel [Sei01a] from a coamoeba point of view. This motivates Conjecture 6.3 in Section 6, which is shown to imply torus-equivariant homological mirror symmetry in Section 7. We illustrate Theorem 1.2 with an example in Section 8. The proof of Theorem 1.4 is given in Section 9.

**Acknowledgment**: We thank Alastair Craw and Akira Ishii for valuable discussions, and the anonymous referee for suggesting several improvements. K. U. is supported by Grant-in-Aid for Young Scientists (No.18840029).

## 2 Dimer models and $A_\infty$-categories

We first recall basic definitions on dimer models:

- A *dimer model* is a bicolored graph $G = (B, W, E)$ on an oriented real 2-torus, which divides the torus into polygons. Here $B$ is the set of black nodes, $W$ is the set of white nodes, and $E$ is the set of edges. No edge is allowed to connect nodes with the same color.

- A *quiver* consists of a finite set $V$ called the set of vertices, another finite set $A$ called the set of arrows, and two maps $s, t : A \to V$ called the source and the target map. The quiver $Q = (V, A, s, t)$ associated with $G$ is defined as the dual graph of $G$, equipped with the orientation so that the white node is always on the right of an arrow; the set $V$ of vertices is the set of faces of $G$, and the set $A$ of arrows can naturally be identified with the set $E$ of edges of $G$.

- A *perfect matching* is a subset $D \subset E$ such that for any $n \in B \sqcup W$, there is a unique edge $e \in D$ adjacent to $n$.

- A perfect matching $D$ is said to be *internal* if there is a choice of a total order $>$ on the set of vertices of the quiver $Q$ such that
$$D = \{a \in A \mid s(a) < t(a)\} \subset A = E. \tag{2.1}$$

- A dimer model is *consistent* if it satisfies the conditions in [IU11, Definition 3.5].

Next we recall the definition of an $A_\infty$-category: An $A_\infty$-category $\mathcal{A}$ consists of

- a set $\mathfrak{Ob}(\mathcal{A})$ of objects,

- for $c_1$, $c_2 \in \mathfrak{Ob}(\mathcal{A})$, a $\mathbb{Z}$-graded vector space $\hom_\mathcal{A}(c_1, c_2)$ called the space of morphisms, and

- operations
$$\mathfrak{m}_l : \hom_\mathcal{A}(c_{l-1}, c_l) \otimes \cdots \otimes \hom_\mathcal{A}(c_0, c_1) \longrightarrow \hom_\mathcal{A}(c_0, c_l)$$
of degree $2 - l$ for $l = 1, 2, \ldots$ and $c_0, \ldots, c_l \in \mathfrak{Ob}(\mathcal{A})$,

satisfying the $A_\infty$-*relations*
$$\sum_{i=0}^{l-1} \sum_{j=i+1}^{l} (-1)^{\deg a_1 + \cdots + \deg a_i - i} \mathfrak{m}_{l+i-j+1}(a_l \otimes \cdots \otimes a_{j+1} \otimes \mathfrak{m}_{j-i}(a_j \otimes \cdots \otimes a_{i+1})$$
$$\otimes a_i \otimes \cdots \otimes a_1) = 0 \tag{2.2}$$

for any positive integer $l$, any sequence $c_0, \ldots, c_l$ of objects of $\mathcal{A}$, and any sequence of morphisms $a_i \in \hom_\mathcal{A}(c_{i-1}, c_i)$ for $i = 1, \ldots, l$.

Let $G = (B, W, E)$ be a dimer model and $Q = (V, A, s, t)$ be the quiver associated with $G$. Then the $A_\infty$-category $\mathcal{A}$ associated with $G$ is defined as follows [FU10, Definition 2.1]:

- The set of objects is the set $V$ of vertices of the quiver.

- For two objects $v$ and $w$ in $\mathcal{A}$, the space of morphisms is given by
$$\hom^i(v, w) = \begin{cases} \mathbb{C} \cdot \mathrm{id}_v & i = 0 \text{ and } v = w, \\ \mathrm{span}\{a \mid a : w \to v\} & i = 1, \\ \mathrm{span}\{a^\vee \mid a : v \to w\} & i = 2, \\ \mathbb{C} \cdot \mathrm{id}_v^\vee & i = 3 \text{ and } v = w, \\ 0 & \text{otherwise.} \end{cases}$$

- Non-zero $A_\infty$-operations are
$$\mathfrak{m}_2(x, \mathrm{id}_v) = \mathfrak{m}_2(\mathrm{id}_w, x) = x$$
for any $x \in \hom(v, w)$,
$$\mathfrak{m}_2(a, a^\vee) = \mathrm{id}_v^\vee$$

and
$$\mathfrak{m}_2(a^\vee, a) = \mathrm{id}_w^\vee$$

for any arrow $a$ from $v$ to $w$,

$$\mathfrak{m}_k(a_1, \ldots, a_k) = a_0$$

for any cycle $(a_0, \ldots, a_k)$ of the quiver going around a white node, and

$$\mathfrak{m}_k(a_1, \ldots, a_k) = -a_0$$

for any cycle $(a_0, \ldots, a_k)$ of the quiver going around a black node.

For an $A_\infty$-category $\mathcal{A}$ and a total order $<$ on the set of objects, the *directed subcategory* $\mathcal{A}^\to$ is the $A_\infty$-category whose set of objects is the same as $\mathcal{A}$ and whose spaces of morphisms are given by

$$\hom_{\mathcal{A}^\to}(X, Y) = \begin{cases} \mathbb{C} \cdot \mathrm{id}_X & X = Y, \\ \hom_\mathcal{A}(X, Y) & X < Y, \\ 0 & \text{otherwise,} \end{cases}$$

with the $A_\infty$-operations inherited from $\mathcal{A}$.

If $(G, D)$ is a pair of a consistent dimer model $G$ and an internal perfect matching $D$ on $G$, then the directed subcategory $\mathcal{A}^\to$ of $\mathcal{A}$ does not depend on the choice of a total order $<$ on vertices of the quiver satisfying (2.1).

## 3 Toric Fano stacks associated with lattice polygons

Let $N = \mathbb{Z}^2$ be a free abelian group of rank two and $M = \mathrm{Hom}(N, \mathbb{Z})$ be the dual group. We write $N_\mathbb{R} = N \otimes \mathbb{R}$, $M_\mathbb{R} = M \otimes \mathbb{R}$ and $\mathbb{T} = N \otimes \mathbb{C}^\times = \mathrm{Spec}\,\mathbb{C}[M]$. Let $\Delta$ be a convex lattice polygon in $N_\mathbb{R}$ containing the origin in its interior, and $\{v_i\}_{i=1}^r$ be the set of vertices of $\Delta$ numbered in accordance with the cyclic order. Let

$$\varphi : \mathbb{Z}^r \to N$$

be the homomorphism of abelian groups sending the $i$-th standard coordinate $e_i \in \mathbb{Z}^r$ to $v_i \in N$ for $i = 1, \ldots, r$. Then the toric Fano stack $X = X_\Delta$ associated with $\Delta$ is the quotient stack

$$X = [(\mathbb{C}^r \setminus \mathcal{SR})/\mathcal{K}]$$

of an open subscheme of $\mathbb{C}^r$ by the natural action of

$$\mathcal{K} = \mathrm{Ker}(\varphi \otimes \mathbb{C}^\times : (\mathbb{C}^\times)^r \to \mathbb{T}),$$

where the Stanley-Reisner locus $\mathcal{SR} \subset \mathbb{C}^r$ is defined by

$$\mathcal{SR} = \{(x_1, \ldots x_r) \in \mathbb{C}^r \mid (x_i, x_j) \neq (0, 0) \text{ if } v_i \text{ and } v_j \text{ are not adjacent}\}$$

if $r > 3$, and consists of the origin if $r = 3$.

Let $\phi: N \to N$ be an endomorphism of rank two. Then the toric Fano stack associated with the polygon $\phi(\Delta)$ is the quotient stack

$$X_{\phi(\Delta)} = [X_\Delta / K]$$

of $X_\Delta$ by the group

$$K = \mathrm{Ker}(\phi \otimes \mathbb{C}^\times : \mathbb{T} \to \mathbb{T}).$$

It follows from [IU, Theorem 7.2] and [FU10, Proposition 3.2] that for a pair $(G, D)$ of a consistent dimer model and an internal perfect matching on it, there is a lattice polygon $\Delta$ containing the origin in its interior such that there is an equivalence

$$D^b \mathcal{A}^\to \cong D^b \mathrm{coh}\, X_\Delta \tag{3.1}$$

of triangulated categories. The lattice polygon $\Delta$ is called the *characteristic polygon* of the pair $(G, D)$, and can be described in a combinatorial way in terms of *height changes* of perfect matchings.

## 4 Fukaya categories

Let $W$ be a regular function on an algebraic torus $(\mathbb{C}^\times)^2 = \mathrm{Spec}\,\mathbb{C}[x^{\pm 1}, y^{\pm 1}]$ equipped with an exact Kähler form

$$\omega = \frac{1}{2\sqrt{-1}} \left( \frac{dx \wedge d\overline{x}}{|x|^2} + \frac{dy \wedge d\overline{y}}{|y|^2} \right).$$

$W$ is a *Lefschetz fibration* if all the critical points are non-degenerate with distinct critical values, and the horizontal lift $\widetilde{\gamma}_x : [0,1] \to X$ of a smooth path $\gamma : [0,1] \to \mathbb{C}$ starting at $x \in p^{-1}(\gamma(0))$ is always defined. The latter condition is satisfied if the Newton polygon $\Delta$ of $W$ contains the origin in its interior, and the leading term of $W$ with respect to any edge of $\Delta$ has no critical point on the torus.

Assume that the origin is a regular value of $W$. A *distinguished basis of vanishing cycles* is a collection $(C_1, \ldots, C_m)$ of embedded circles in $W^{-1}(0)$ which collapse to critical points by parallel transport along a *distinguished set of vanishing paths*, cf. [Sei08, Section 16]. We assume that vanishing cycles intersect each other transversely.

A *relative grading* of $W$ is a nowhere-vanishing smooth section of the holomorphic line bundle $\Lambda^{\mathrm{top}}(T^*(\mathbb{C}^\times)^2)^{\otimes 2} \otimes W^*(T^*\mathbb{C})^{\otimes (-2)}$, which we choose as $(d\log x \wedge d\log y)^{\otimes 2} \otimes (dW)^{\otimes (-2)}$. It induces a section $\eta$ of $\Lambda^{\mathrm{top}}(T^*M)^{\otimes 2}$ on the fiber $M = W^{-1}(0)$, which gives a map

$$\begin{array}{rccc}
\det^2_\eta : & \mathcal{L}ag_M & \to & \mathbb{C}^\times / \mathbb{R}^{>0} \cong S^1 \\
& \cup & & \cup \\
& \mathrm{span}\{e_1, \ldots, e_n\} & \mapsto & [\eta((e_1 \wedge \cdots \wedge e_n)^{\otimes 2})]
\end{array}$$

from the Lagrangian Grassmannian bundle $\mathcal{L}ag_M$ on $M$. A Lagrangian submanifold $L \subset M$ naturally gives a section

$$\begin{array}{rccc}
s_L : & L & \to & \mathcal{L}ag_M|_L \\
& \cup & & \cup \\
& x & \mapsto & T_x L.
\end{array}$$

A grading of $L$ is a lift $\widetilde{\phi}_L : L \to \mathbb{R}$ of the composition

$$\phi_L = \det_\eta^2 \circ s_L : L \to S^1$$

to the universal cover $\mathbb{R} \to S^1$.

Given a pair $(L_1, L_2)$ of graded Lagrangian submanifolds, one can define the *Maslov index* $\mu(x; L_1, L_2)$ for each intersection point $x \in L_1 \cap L_2$. Since $\dim_\mathbb{C} M = 1$, it is given by the round-up

$$\mu(x; L_1, L_2) = \lfloor \widetilde{\phi}_{L_2}(x) - \widetilde{\phi}_{L_1}(x) \rfloor$$

of the difference of the phase functions at $x$.

The Fukaya category $\mathfrak{Fuk}\, W$ of $W$ is a directed $A_\infty$-category whose set of object is $(C_i)_{i=1}^m$ and whose spaces of morphisms are given by

$$\hom(C_i, C_j) = \begin{cases} \mathrm{id}_{C_i} & i = j, \\ \bigoplus_{p \in C_i \cap C_j} \mathbb{C} \cdot p & i < j, \\ 0 & \text{otherwise.} \end{cases}$$

The $A_\infty$-operation $\mathfrak{m}_k$ involving the identity morphism is the obvious one for $k=2$ and zero for $k \neq 2$. Non-trivial $A_\infty$-operations are given by

$$\mathfrak{m}_k(p_k, \ldots, p_1) = \sum_{p_0 \in \cap C_{i_k} \cap C_{i_0}} (-1)^\dagger p_0, \qquad (4.1)$$

where $p_\ell \in C_{i_{\ell-1}} \cap C_{i_\ell}$ for $\ell = 1, \ldots, k$ and the sum is over the '$(k+1)$-gons' whose $\ell$-th vertex is the point $p_\ell$ for $\ell = 0, \ldots, k$ and whose edge between $p_\ell$ and $p_{\ell+1}$ lies on $C_{i_\ell}$. The grading of $C_{i_\ell}$ defines an orientation of $C_{i_\ell}$, and let $\xi_\ell$, $\ell = 0, \ldots, k$ be the unit tangent vector of $C_{i_\ell}$ at $p_\ell$ along the orientation. We also choose a branch point on each vanishing cycle $C_{i_\ell}$, which comes from the choice of the non-trivial spin structure. In this paper, we only deal with the case where $\mu(p_\ell; C_{i_{\ell-1}}, C_{i_\ell}) = 1$ for $\ell = 1, \ldots, k$ and $\mu(p_0; C_{i_0}, C_{i_k}) = 2$. In this case, the sign rule of Seidel [Sei03, Section (9e)] states that † in (4.1) is the sum of (i) the number of $1 \leq \ell \leq k$ such that $\xi_\ell$ points away from the $(k+1)$-gon, and (ii) the the number of branch points on the $(k+1)$-gon coming from the spin structures.

## 5 Coamoeba for the mirror of $\mathbb{P}^2$

The mirror of $\mathbb{P}^2$ is given by the Laurent polynomial

$$W(x, y) = x + y + \frac{1}{xy}.$$

The critical points of $W$ are given by

$$(x, y) = (1, 1), (\omega, \omega), (\omega^2, \omega^2),$$

where $\omega = \exp(2\pi \sqrt{-1}/3)$ is a primitive cubic root of unity. The corresponding critical values are $3, 3\omega$ and $3\omega^2$. Let $(c_i)_{i=1}^3$ be the distinguished set of vanishing paths obtained as the straight line segments from the origin to the critical values of $W$ as shown in Figure

5.1. The fiber $W^{-1}(0)$ can be realized as a branched double cover of the $x$-plane by the projection
$$\begin{array}{rcl} \pi: W^{-1}(0) & \to & \mathbb{C}^\times \\ \cup & & \cup \\ (x, y) & \mapsto & x. \end{array}$$

The branch points of $\pi$ are given by three dots in Figure 5.2, together with $x = 0$ which does not belong to $\mathbb{C}^\times$. The fiber $W^{-1}(0)$ can be compactified to an elliptic curve by adding one point over $x = 0$ and two points over $x = \infty$.

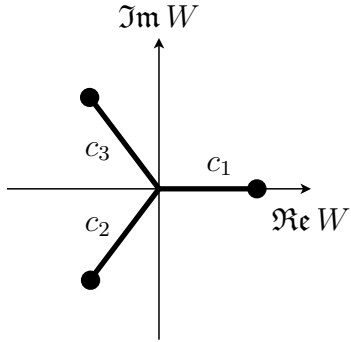
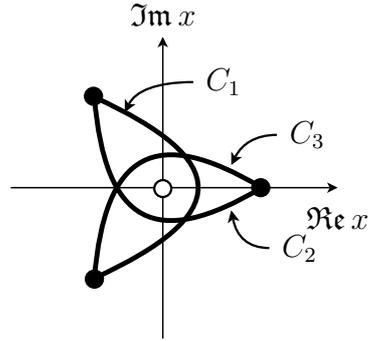

Figure 5.1: A path on the $W$-plane

Figure 5.2: The trajectories of the branch points

Recall that the *coamoeba* of a subvariety of the torus $(\mathbb{C}^\times)^2$ is defined by Passare and Tsikh as its image by the argument map
$$\begin{array}{rcl} \mathrm{Arg}: (\mathbb{C}^\times)^2 & \to & \mathbb{R}^2/\mathbb{Z}^2 \\ \cup & & \cup \\ (x, y) & \mapsto & \dfrac{1}{2\pi}(\arg x, \arg y). \end{array}$$

It follows from [UY11, Theorem 7.1] that the coamoeba of $W^{-1}(0)$ is the union of interiors and vertices of six triangles in Figure 5.3. The inverse image of the set of vertices divide $W^{-1}(0)$ into six triangles as in Figure 5.4. Figures 5.5 and 5.6 show these triangles on two sheets.

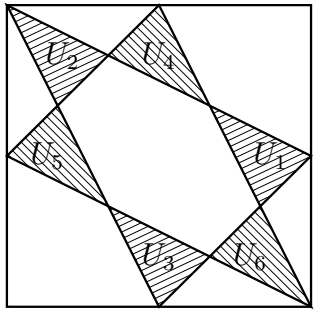
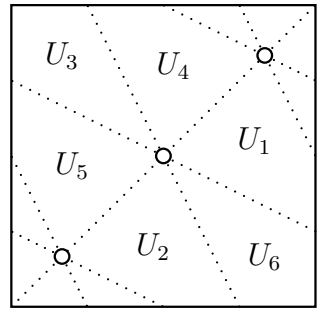

Figure 5.3: The coamoeba of $W^{-1}(0)$

Figure 5.4: The fiber $W^{-1}(0)$

If one deforms $W^{-1}(0)$ to $W^{-1}(c_i(t))$ for $i = 1, 2, 3$ and $t \in [0, 1]$, the branch points move as in Figure 5.2. Here, $C_1$, $C_2$, and $C_3$ are trajectories of the branch points along $c_1$, $c_2$ and $c_3$ respectively. These trajectories are the images of the corresponding vanishing cycles by $\pi$ up to homotopy. This shows that the vanishing cycles on $W^{-1}(0)$ are as in Figure 5.7. We choose the gradings and the spin structures on these vanishing cycles in such a way that the arrows show the orientations of vanishing cycles coming from gradings and the dots show the branch points for the non-trivial spin structures. The images of vanishing cycles under the argument map are shown in Figure 5.8.

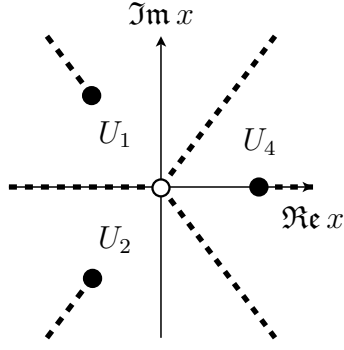

Figure 5.5: The first sheet

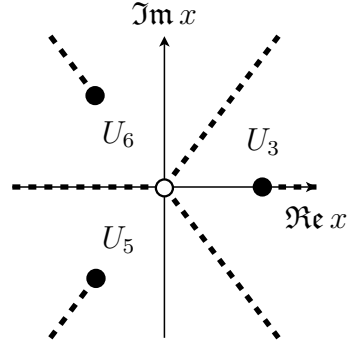

Figure 5.6: The second sheet

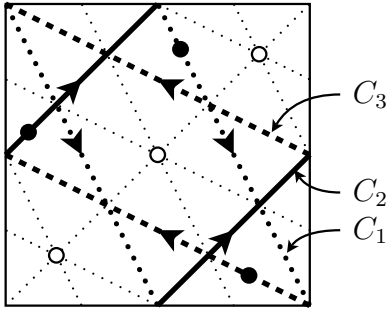

Figure 5.7: Vanishing cycles on $W^{-1}(0)$

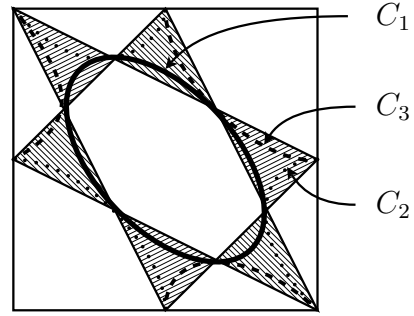

Figure 5.8: Vanishing cycles on the coamoeba

One can see from Figure 5.7 that six triangles in Figure 5.9 are the only polygons on $W^{-1}(0)$ bounded by $\bigcup_i C_i$. By contracting these triangles, one obtains the graph on $W^{-1}(0)$ shown in Figure 5.10. Non-trivial $A_\infty$-operations on $\mathfrak{Fuk}\, W$ are in one-to-one correspondence with these triangles, and the colors of the nodes come from the signs of these $A_\infty$-operations determined by Seidel's rule recalled in Section 4; the sign is positive for a white node and negative for a black node.

The argument projection $G$ of the graph in Figure 5.10 is shown in Figure 5.11, which gives a consistent dimer model. There is a natural bijection between the set of faces of $G$ and the distinguished basis $(C_i)_{i=1}^3$ of vanishing cycles. Moreover, intersections of two vanishing cycles correspond to common edges of two faces under this bijection.

The order on the distinguished basis of vanishing cycles defines the internal perfect matching in Figure 5.12 by (2.1). The gradings on $C_i$ are chosen so that the Maslov index $\mu(p; C_i, C_j)$ for $p \in C_i \cap C_j$ and $i < j$ is one if the edge corresponding to $p$ is not contained in $D$, and two otherwise.

This suffices to show the equivalence

$$\mathfrak{Fuk}\, W \cong \mathcal{A}^{\rightarrow}$$

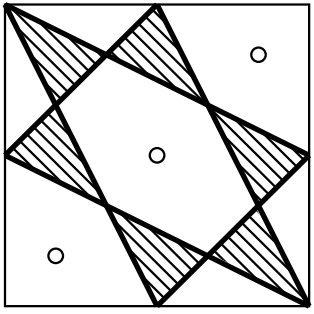
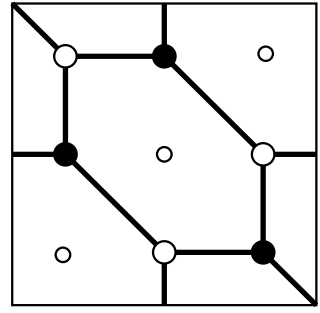

Figure 5.9: Six triangles on $W^{-1}(0)$      Figure 5.10: A graph on $W^{-1}(0)$

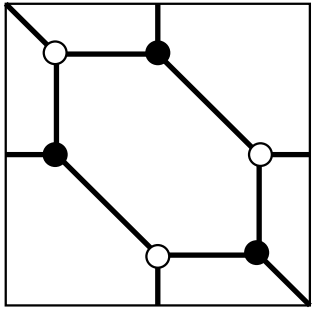
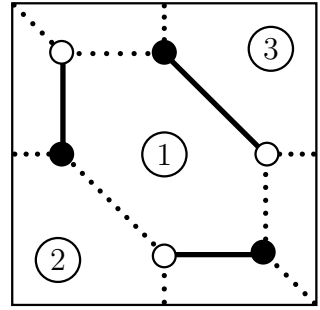

Figure 5.11: A graph on $T$      Figure 5.12: The perfect matching

between the Fukaya category of $W$ and the directed $A_\infty$-category associated with the pair $(G, D)$ defined in Section 2: Faces and edges of $G$ correspond to objects and morphisms of both $\mathfrak{Fuk}\, W$ and $\mathcal{A}^\rightarrow$. Non-trivial $A_\infty$-operations in $\mathfrak{Fuk}\, W$ come from triangles bounded by $\bigcup_i C_i$, which are in one-to-one correspondence with nodes of $G$. Non-trivial $A_\infty$-operations in $\mathcal{A}^\rightarrow$ also come from nodes of $G$ by definition, and exactly matches the non-trivial $A_\infty$-operations in $\mathfrak{Fuk}\, W$.

The characteristic polygon of the pair $(G, D)$ coincides with the Newton polygon $\Delta$ of $W$, and the toric Fano stack associated with $\Delta$ is the projective plane $\mathbb{P}^2$. Combined with the equivalence
$$D^b \mathcal{A}^\rightarrow \cong D^b \operatorname{coh} \mathbb{P}^2$$
in (3.1), homological mirror symmetry for $\mathbb{P}^2$ is proved.

# 6 Dimer models as spines of coamoebas

The discussion in Section 5 motivates the following definition:

**Definition 6.1.** Let
$$W : \mathbb{T}^\vee \to \mathbb{C}$$
be a regular map on an algebraic torus $\mathbb{T}^\vee = M \otimes \mathbb{C}^\times = \operatorname{Spec} \mathbb{C}[N]$, $G$ be a consistent dimer model on $T^\vee = M_\mathbb{R}/M$, and $D$ be an internal perfect matching on $G$. Then the pair $(G, D)$ is said to be *associated with* $W$ if there is a distinguished basis $(C_i)_{i=1}^m$ of vanishing cycles in $W^{-1}(0)$ satisfying the following conditions:

- For any point $p \in W^{-1}(0) \setminus \bigcup_i C_i$, there is at most one polygon bounded by $\bigcup_i C_i$ and passing through $p$.

- By contracting these polygons, one obtains a bipartite graph $Y$ on $W^{-1}(0)$.

- The restriction of the argument map $\mathrm{Arg} : \mathbb{T}^\vee \to T^\vee$ to the graph $Y$ is injective.

- The image $\mathrm{Arg}(Y)$ is a consistent dimer model $G$ on $T^\vee$ with respect to a suitable choice of colors on the nodes.

- There is a natural bijection between the set of faces of $G$ and vanishing cycles.

- Intersections of two vanishing cycles corresponds to common edges of two faces under the above bijection.

- The order on the distinguished basis of vanishing cycles determines an internal perfect matching by (2.1).

- With a suitable choice of gradings on $C_i$, the Maslov index $\mu(p; C_i, C_j)$ of an intersection $p \in C_i \cap C_j$ for $i < j$ corresponding to an edge $e$ of $G$ is given by

$$\mu(p; C_i, C_j) = \begin{cases} 1 & e \notin D, \\ 2 & e \in D. \end{cases}$$

- For a suitable choice of branch points for the non-trivial spin structures on $C_i$, the contribution of a polygon in $W^{-1}(0)$ corresponding to a node of $G$ to the $A_\infty$-operations in $\mathfrak{Fuk}\, W$ is $+1$ or $-1$ depending on the color of the node.

- The characteristic polygon of the pair $(G, D)$ coincides with the Newton polygon of $W$.

The conditions in Definition 6.1 ensure the following:

**Proposition 6.2.** *Let $(G, D)$ be a pair of a consistent dimer model $G$ and an internal perfect matching $D$ associated with $W$. Then one has an equivalence*

$$\mathcal{A}^\to \cong \mathfrak{Fuk}\, W$$

*of $A_\infty$-categories.*

*Proof.* This is obtained by comparing the definition of $\mathcal{A}^\to$ given in Section 2 and the definition of $\mathfrak{Fuk}\, W$ in Section 4 using the conditions in Definition 6.1. The set of objects of $\mathfrak{Fuk}\, W$ is a distinguished basis of vanishing cycles, which corresponds to the set of faces of $G$ by Definition 6.1. The set of faces of $G$ in turn corresponds to the set of objects of $\mathcal{A}^\to$ by the definition of $\mathcal{A}^\to$. The spaces of morphisms in $\mathfrak{Fuk}\, W$ are spanned by intersection points of vanishing cycles, which correspond to edges of $G$. The edges of $G$ in turn span the spaces of morphisms in $\mathcal{A}^\to$. The nodes of $G$ are in bijective correspondence with polygons in $W^{-1}(0)$ bounded by $\bigcup_i C_i$, which give non-trivial $A_\infty$-operations in $\mathfrak{Fuk}\, W$. The coefficient for a non-trivial $A_\infty$-operation is $+1$ or $-1$ depending on the color of the corresponding node. This exactly matches the $A_\infty$-operations on $\mathcal{A}^\to$ defined in Section 2, and Proposition 6.2 is proved. □

Now we can state the main conjecture in this paper:

**Conjecture 6.3.** *Let $\Delta$ be a convex lattice polygon in $N_{\mathbb{R}}$ which contains the origin in its interior. Then there exist*

- *a Laurent polynomial $W \in \mathbb{C}[N]$ whose Newton polygon coincides with $\Delta$, and*
- *a pair $(G, D)$ of a consistent dimer model $G$ on $T^{\vee}$ and an internal perfect matching $D$ on $G$ associated with $W$.*

The equivalence (3.1) and Proposition 6.2 show the following:

**Proposition 6.4.** *Conjecture 6.3 for a polygon $\Delta$ implies homological mirror symmetry*

$$D^b \operatorname{coh} X_{\Delta} \cong D^b \mathfrak{Fuk} W$$

*for the toric Fano stack $X_{\Delta}$ associated with $\Delta$.*

# 7 Torus-equivariant homological mirror symmetry

In this section, we discuss the relation between dimer models and torus-equivariant homological mirror symmetry for two-dimensional toric Fano stacks:

**Proposition 7.1.** *Let $\Delta$ be a convex lattice polygon in $N_{\mathbb{R}}$ containing the origin in its interior, and $\phi : N \to N$ be an integral linear transformation of rank two. If Conjecture 6.3 holds for $\Delta$, then it also holds for the polygon $\phi(\Delta)$.*

*Proof.* Let $W$ be a Laurent polynomial whose Newton polygon coincides with $\Delta$ and $(G, D)$ be a pair associated with $W$. For an integral linear transformation $\phi$, we have to show the existence of another Laurent polynomial $W'$ and a pair $(G', D')$ such that the Newton polygon of $W'$ is $\phi(\Delta)$ and the pair $(G', D')$ is associated with $W'$.

Let $\psi : M \to M$ be the transpose of $\phi$ and

$$W' = W \circ (\psi \otimes \mathbb{C}^{\times}) : \mathbb{T}^{\vee} \to \mathbb{C}$$

be the pull-back of $W$ by the the $\det(\psi)$-fold covering map $\psi \otimes \mathbb{C}^{\times} : \mathbb{T}^{\vee} \to \mathbb{T}^{\vee}$. Similarly, define the pair $(G', D')$ as the pull-back of the pair $(G, D)$ on $T^{\vee}$ by the $\det(\psi)$-fold covering map $\psi \otimes (\mathbb{R}/\mathbb{Z}) : T^{\vee} \to T^{\vee}$. The numbers of nodes and edges of $G'$ are $\det(\psi)$ times the numbers of nodes and edges of $G$. Since $\psi \otimes \mathbb{C}^{\times} : \mathbb{T}^{\vee} \to \mathbb{T}^{\vee}$ is an unramified covering of $\mathbb{T}^{\vee}$, the set of critical values of $W'$ coincides with those of $W$, and a distinguished basis of vanishing cycles of $W'$ is obtained as the pull-back of a distinguished basis of vanishing cycles of $W$. The fact that there are as many as $\det(\psi)$ vanishing cycles for a given critical value does not cause any problem, since these $\det(\psi)$ vanishing cycles are mutually disjoint. The commutativity of the diagram

$$\begin{array}{ccc} \mathbb{T}^{\vee} & \xrightarrow{\psi \otimes \mathbb{C}^{\times}} & \mathbb{T}^{\vee} \\ \scriptsize{\operatorname{Arg}} \downarrow & & \scriptsize{\operatorname{Arg}} \downarrow \\ T^{\vee} & \xrightarrow{\psi \otimes (\mathbb{R}/\mathbb{Z})} & T^{\vee} \end{array}$$

and the fact that everything is pulled-back from the column on the right by the unramified covering maps $\psi \otimes \mathbb{C}^{\times}$ and $\psi \otimes (\mathbb{R}/\mathbb{Z})$ show that the pair $(G', D')$ is associated with $W'$ in the sense of Definition 6.1, and Proposition 7.1 is proved. □

Propositions 6.4 and 7.1 imply Theorem 1.3. Theorem 1.2 is obtained in a similar way from Propositions 7.2 and 7.3 below:

**Proposition 7.2.** *Let $(G, D)$ be a pair of a consistent dimer model and an internal perfect matching, $\mathcal{A}^\to$ be the directed $A_\infty$-category associated with $(G, D)$, and $X$ be the toric Fano stack associated with the characteristic polygon of $(G, D)$. Let further $\widetilde{\mathcal{A}}^\to$ be the directed $A_\infty$-category associated with the pull-back $(\widetilde{G}, \widetilde{D})$ of $(G, D)$ by the universal cover $M_\mathbb{R} \to T^\vee = M_\mathbb{R}/M$. Then one as an equivalence*
$$D^b \widetilde{\mathcal{A}}^\to \cong D^b \operatorname{coh}^\mathbb{T} X$$
*of triangulated categories.*

Although the pair $(\widetilde{G}, \widetilde{D})$ is infinite, the definition of the directed $A_\infty$-category associated with a pair given in Section 2 makes sense also in this case, and the proof of the equivalence (3.1) in [FU10, Proposition 3.2] carries over verbatim to the proof of Proposition 7.2.

**Proposition 7.3.** *Let $(G, D)$ be a pair of a consistent dimer model and an internal perfect matching associated with a Laurent polynomial $W$. Then the pull-back $(\widetilde{G}, \widetilde{D})$ of $(G, D)$ by the universal cover $M_\mathbb{R} \to T^\vee$ is associated with $\widetilde{W} = W \circ \exp : \mathbb{C}^2 \to \mathbb{C}$, so that one has an equivalence*
$$D^b \widetilde{\mathcal{A}}^\to \cong D^b \mathfrak{Fuk}\, \widetilde{W}$$
*of triangulated categories.*

The proof of Proposition 7.3 is obtained from the proof of Proposition 7.1 by replacing the $\det(\psi)$-fold coverings $\psi \otimes \mathbb{C}^\times : \mathbb{T}^\vee \to \mathbb{T}^\vee$ and $\psi \otimes (\mathbb{R}/\mathbb{Z}) : T^\vee \to T^\vee$ with the universal coverings $\exp : M_\mathbb{C} \to \mathbb{T}^\vee$ and $M_\mathbb{R} \to T^\vee$ respectively.

# 8 An example of a quotient stack of $\mathbb{P}^2$

In this section, we illustrate Theorem 1.3 with an example of a quotient stack of $\mathbb{P}^2$ by a cyclic group of order three. Let $\Delta$ be the convex lattice polygon in Figure 8.1 corresponding to $\mathbb{P}^2$ and $\phi : N \to N$ be a linear map represented by the matrix $\begin{pmatrix} 2 & -1 \\ -1 & 2 \end{pmatrix}$. Then $\phi(\Delta)$ is the convex hull of $(2, -1)$, $(-1, 2)$ and $(-1, -1)$ shown in Figure 8.2.

One has
$$K = \operatorname{Ker} \phi \otimes \mathbb{C}^\times = \langle (\omega, \omega^2) \rangle \cong \mathbb{Z}/3\mathbb{Z}$$
in this case, where $\omega = \exp(2\pi\sqrt{-1}/3)$ is a primitive cubic root of unity. We choose

$$W_1(x, y) = x + y + \frac{1}{xy},$$
$$W_2(x, y) = \frac{x^2}{y} + \frac{y^2}{x} + \frac{1}{xy},$$

so that the Newton polygons of $W_1$ and $W_2$ are given by $\Delta$ and $\phi(\Delta)$ respectively, and one has
$$W_2 = W_1 \circ (\psi \otimes \mathbb{C}^\times)$$

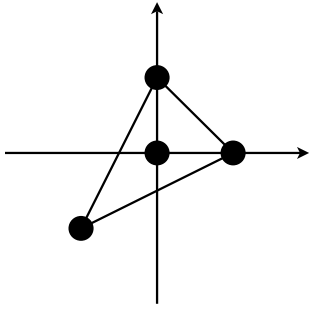

Figure 8.1: The lattice polygon $\Delta$

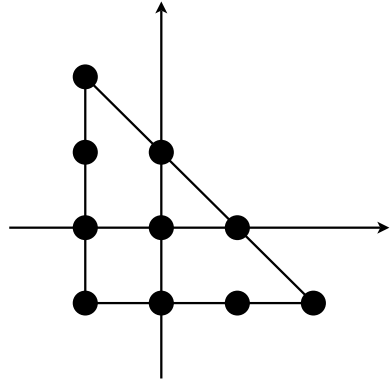

Figure 8.2: The lattice polygon $\phi(\Delta)$

where $\psi : M \to M$ is the transpose of $\phi$. Figure 8.3 show the pull-back of the dimer model in Figure 5.11 associated with $W_1$ by the projection $\pi : M_\mathbb{R} \to T^\vee$, and the dimer model associated with $W_2$ is obtained as its quotient by the lattice $\psi(M)$ generated by $(2, -1)$ and $(-1, 2)$. The resulting dimer model is shown in Figure 8.4, which is a covering of the former with the group of characters of $K$ as the group of deck transformations.

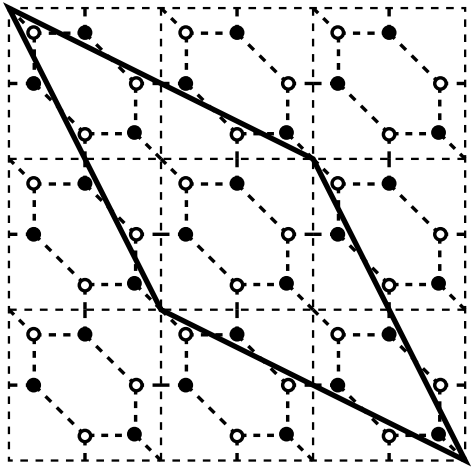

Figure 8.3: The pull-back by $\pi$ of the dimer model associated with $W_1$ and a fundamental region of the action by the lattice $\psi(M)$

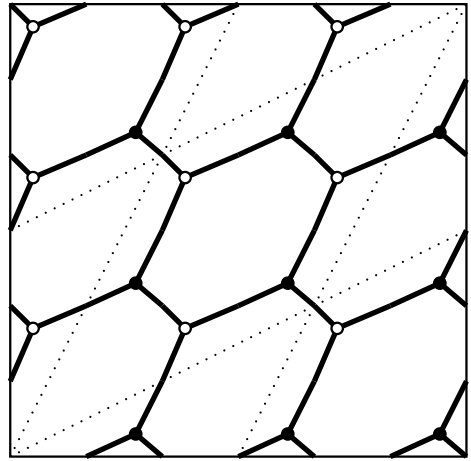

Figure 8.4: A dimer model associated with $W_2$

Figure 8.5 shows the correspondence between fundamental regions of the $M$-action and characters of the torus $\mathbb{T}$. The operation of taking the quotient with respect to the action of $\psi(M)$ corresponds to the restriction of characters to $K \subset \mathbb{T}$, and the resulting characters of $K$ are shown in Figure 8.6.

By comparing Figures 5.12, 8.3 and 8.6, one can label faces of the dimer model in Figure 8.4 as in Figure 8.7; the face labeled as $(ij)$ corresponds to the $j$-th lift to $W_2^{-1}(0)$ of the vanishing cycle $C_i$ of $W_1$.

The fiber $W_2^{-1}(0)$ of $W_2$ can be obtained by assigning a triangle to each node in Figure 8.8 and gluing them together as in Figure 8.9. Vanishing cycles on the resulting elliptic curve minus nine points are shown in Figure 8.10.

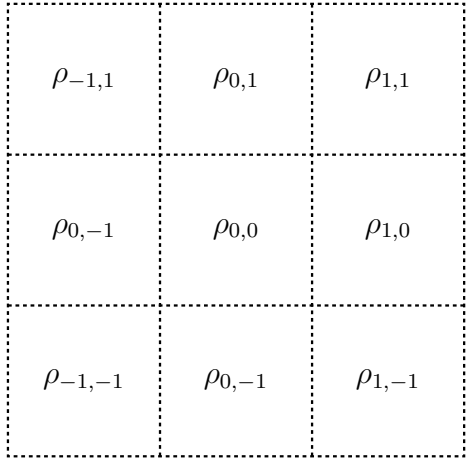

Figure 8.5: Characters of $\mathbb{T}$

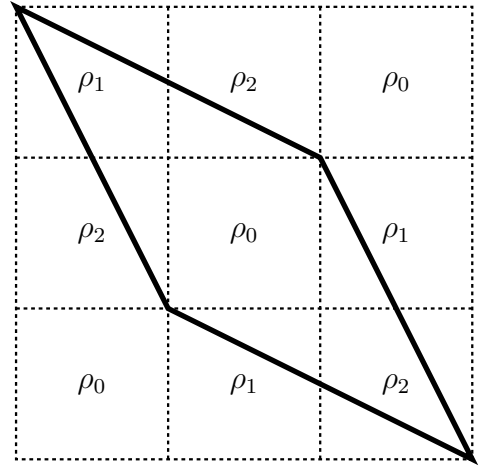

Figure 8.6: Characters of $K$

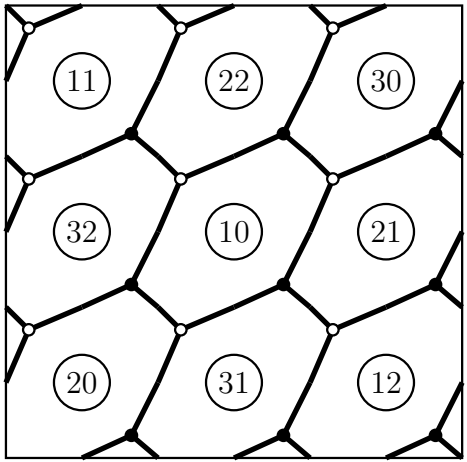

Figure 8.7: Labels on the faces

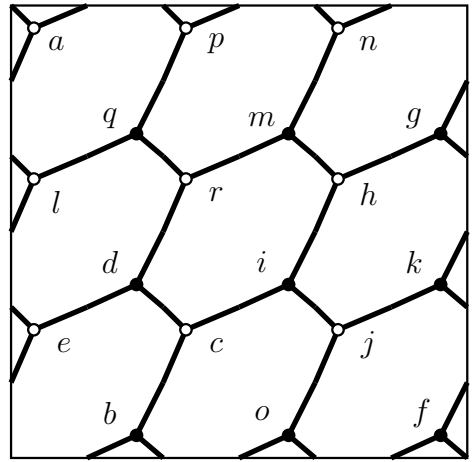

Figure 8.8: Labels on the nodes

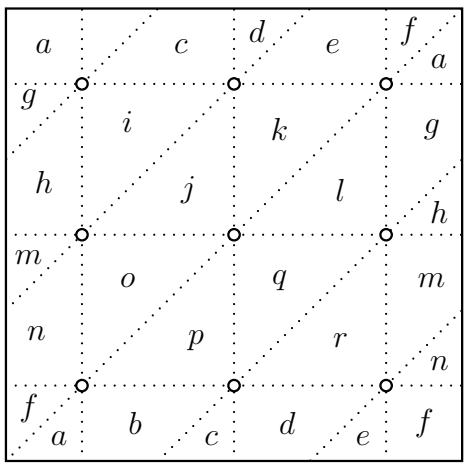

Figure 8.9: The fiber $W_2^{-1}(0)$

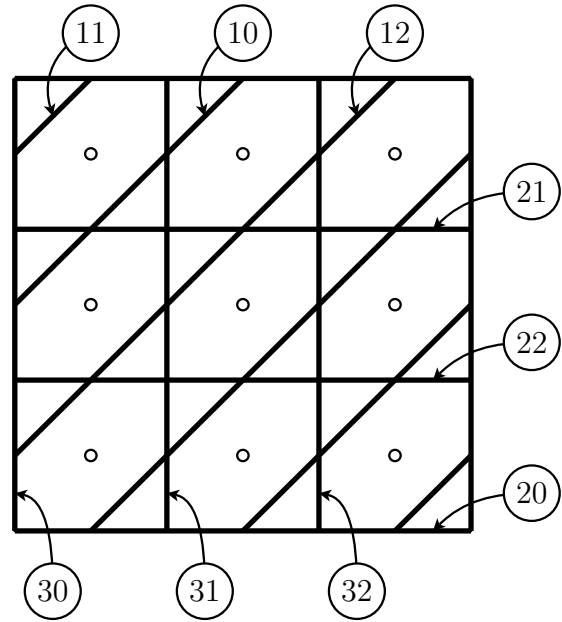

Figure 8.10: Vanishing cycles on $W_2^{-1}(0)$

Let $(E_1, E_2, E_3) = (\Omega^2_{\mathbb{P}^2}(2)[2], \Omega^1_{\mathbb{P}^2}(1)[1], \mathcal{O}_{\mathbb{P}^2})$ be the full exceptional collection in $D^b \operatorname{coh} \mathbb{P}^2$ corresponding to the distinguished basis $(C_1, C_2, C_3)$ of vanishing cycles of $W_1$ under homological mirror symmetry. Then the object $E_i \otimes \rho_j$ in $D^b \operatorname{coh}[\mathbb{P}^2/K] \cong D^b \operatorname{coh}^K \mathbb{P}^2$ correspond to the face labeled as $\boxed{ij}$ in Figure 8.7, which in turn corresponds to the vanishing cycle $\boxed{ij}$ on $W_2^{-1}(0)$ shown in Figure 8.10. Non-trivial $A_\infty$-operations in both $\mathfrak{Fuk}\, W_2$ and $D^b \operatorname{coh}[\mathbb{P}^2/K]$ are in one-to-one correspondence with nodes of the dimer model shown in Figure 8.8. In this way, the pair $(G, D)$ associated with $W_1$ encodes homological mirror symmetry for any toric orbifolds of $\mathbb{P}^2$.

# 9 Toric del Pezzo surfaces

We prove Theorem 1.4 in this section. We first discuss the case when $X$ is $\mathbb{P}^2$ blown-up at one point. The corresponding lattice polygon $\Delta$ is the convex hull of

$$v_1 = (1, 0), v_2 = (0, 1), v_3 = (-1, 0), \text{ and } v_4 = (-1, -1).$$

Take the Laurent polynomial

$$W(x, y) = x + y - \frac{1}{x} + \frac{1}{xy},$$

whose Newton polygon coincides with $\Delta$. Let $(c_i)_{i=1}^4$ be the distinguished set of vanishing paths, defined as the straight line segments from the origin to the critical values as in Figure 9.1. The second projection will be denoted by

$$\begin{array}{rcl} \pi \colon W^{-1}(0) & \to & \mathbb{C}^\times \\ \cup & & \cup \\ (x, y) & \mapsto & y. \end{array}$$

The fiber $\pi^{-1}(y)$ consists of two points for $y \in \mathbb{C}^\times \setminus \{1\}$, whereas $\pi^{-1}(1)$ consists of one point (the other point goes to $x = 0$). Figure 9.2 shows the image by $\pi$ of the distinguished basis $(C_i)_{i=1}^4$ of vanishing cycles along the paths $(c_i)_{i=1}^4$. The black dots are the branch points, the white dots are $y = 0$ and $y = 1$, the solid lines are the images of the vanishing cycles and the dotted lines are cuts introduced artificially to divide $W^{-1}(0)$ into two sheets as shown in Figure 9.5 and Figure 9.6.

Figure 9.3 shows a schematic picture of the coamoeba. To study the image of the vanishing cycles by the argument map, we cut the coamoeba into pieces along the bold lines in Figure 9.4. These lines look as in Figure 9.5 and Figure 9.6 on the two sheets. They cut $W^{-1}(0)$ into the union of two quadrilaterals and four triangles, glued along ten edges. By gluing these pieces, one obtains an elliptic curve minus four points shown in Figure 9.7. One can see that the vanishing cycles on $W^{-1}(0)$ look as in Figure 9.8 using Figures 9.2, 9.5, 9.6, and 9.7. The black dots in Figure 9.8 are branch points for the non-trivial spin structures on vanishing cycles, and the arrows show the orientations of the vanishing cycles coming from a choice of gradings on vanishing cycles. Figure 9.9 shows the argument projections of vanishing cycles. One can see in Figure 9.8 that there are two quadrangles and four triangles bounded by vanishing cycles, which give rise to non-trivial

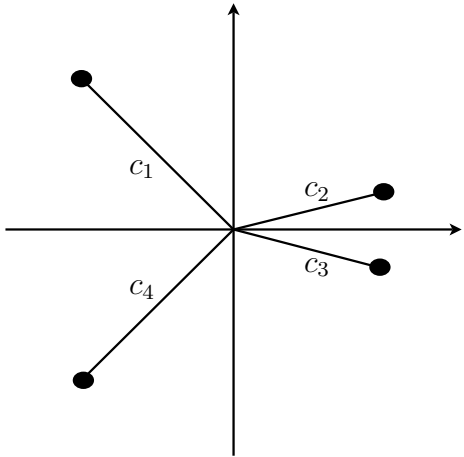

Figure 9.1: A distinguished basis of vanishing paths

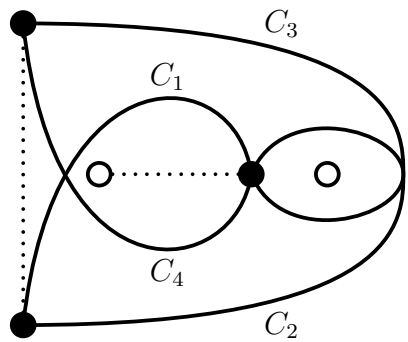

Figure 9.2: Vanishing cycles on the $y$-plane

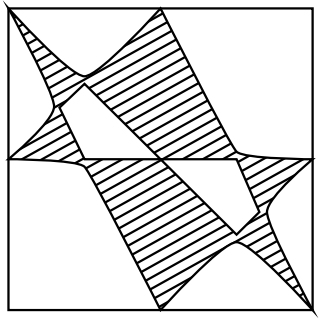

Figure 9.3: The coamoeba of $W^{-1}(0)$

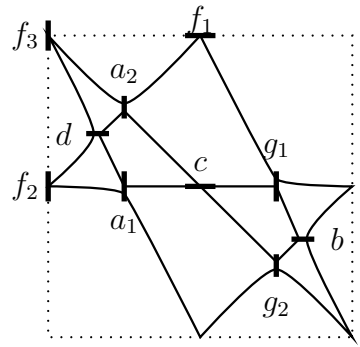

Figure 9.4: Cutting the coamoeba into pieces

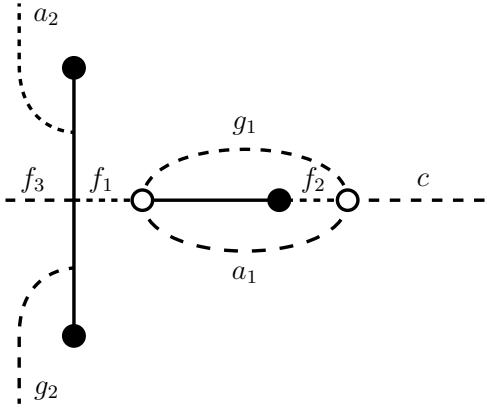

Figure 9.5: The first sheet of $W^{-1}(0)$

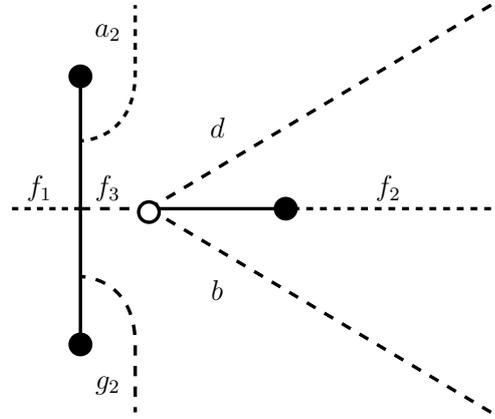

Figure 9.6: The second sheet of $W^{-1}(0)$

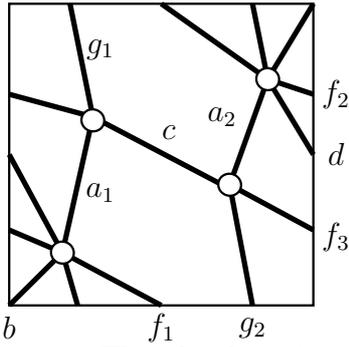

Figure 9.7: The glued surface

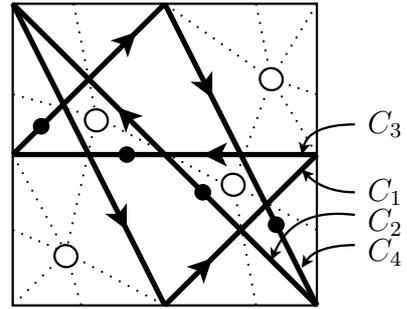

Figure 9.8: The vanishing cycles

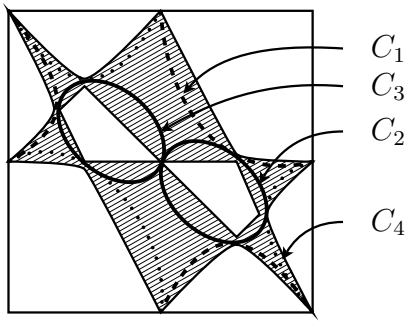

Figure 9.9: Vanishing cycles on the coamoeba

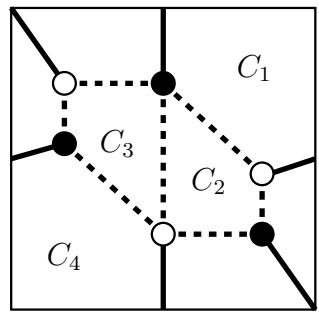

Figure 9.10: The dimer model and the perfect matching

$A_\infty$-operations in the Fukaya category. By contracting these polygons, one obtains a graph on $W^{-1}(0)$, which projects to a dimer model by the argument map. Figure 9.10 shows the pair $(G, D)$ of the resulting dimer model $G$ and a perfect matching $D$ on it, where the solid lines and dotted lines show edges of $G$ belonging to $D$ and $E \setminus D$ respectively. The colors of the nodes come from the signs in the corresponding $A_\infty$-operations. This shows that Conjecture 6.3 holds in this case.

The case of $\mathbb{P}^2$ blown-up at two or three points is completely parallel to the case of $\mathbb{P}^2$ blown-up at one point above. In the case of $\mathbb{P}^2$ blown-up at two points, we choose the Laurent polynomial

$$W(x,y) = (0.8 + 0.6\sqrt{-1})x + (0.8 - 0.6\sqrt{-1})y + \frac{1}{x} + \frac{1}{xy} + \frac{1}{y},$$

which has five non-degenerate critical values. We choose straight line segments from the origin to critical values shown in Figure 9.11 as a distinguished set of vanishing paths. Figure 9.12 shows the images of the corresponding vanishing cycles by the $y$-projection. Vanishing cycles on $W^{-1}(0)$ are as in Figure 9.13. One can see that there are one pentagon, two quadrangles, and three triangles bounded by vanishing cycles, which contribute to the $A_\infty$-operations in the Fukaya category. By contracting these polygons and taking the argument projection of the resulting graph, one obtains the pair $(G, D)$ of a dimer model $G$ and a perfect matching $D$ shown in Figure 9.14.

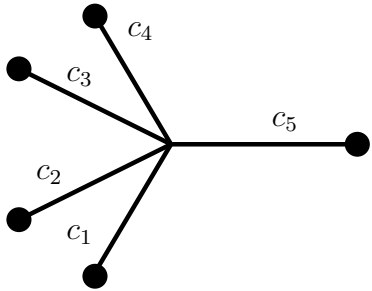

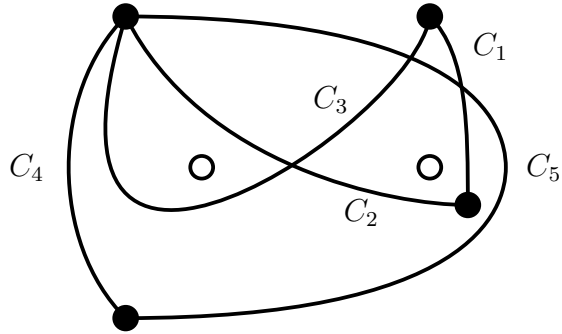

Figure 9.11: A distinguished set of vanishing paths

Figure 9.12: Vanishing cycles on the $y$-plane

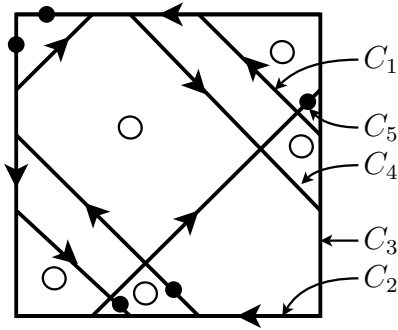

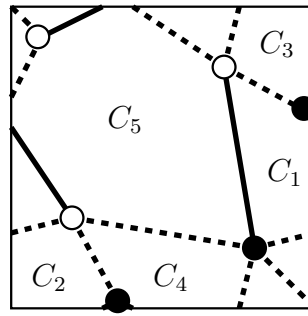

Figure 9.13: Vanishing cycles on $W^{-1}(0)$

Figure 9.14: The dimer model and the perfect matching

In the case of $\mathbb{P}^2$ blown-up at three points, we choose the Laurent polynomial

$$W_6(x,y) = -x - \sqrt{-1}xy - \sqrt{-1}y - \frac{\sqrt{-1}}{x} - \frac{\sqrt{-1}}{xy} + \frac{1}{y},$$

which has six non-degenerate critical values. We choose straight line segments from the origin to critical values shown in Figure 9.15 as a distinguished set of vanishing paths. Figure 9.16 shows the images of the corresponding vanishing cycles by the $y$-projection. Vanishing cycles on $W^{-1}(0)$ are as in Figure 9.17. One can see that there are one hexagon, two quadrangles, and three triangles bounded by vanishing cycles, which contribute to the $A_\infty$-operations in the Fukaya category. By contracting these polygons and taking the argument projection of the resulting graph, one obtains the pair $(G, D)$ of a dimer model $G$ and a perfect matching $D$ shown in Figure 9.18.

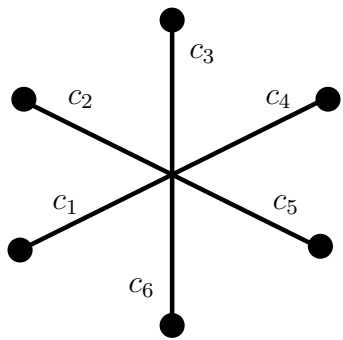

Figure 9.15: A distinguished set of vanishing paths

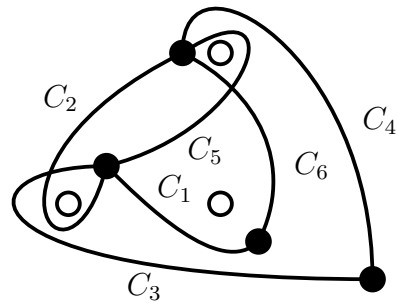

Figure 9.16: Vanishing cycles on the $y$-plane

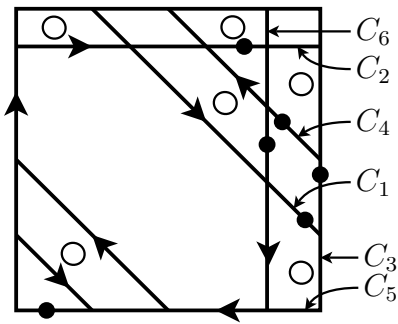

Figure 9.17: Vanishing cycles on $W^{-1}(0)$

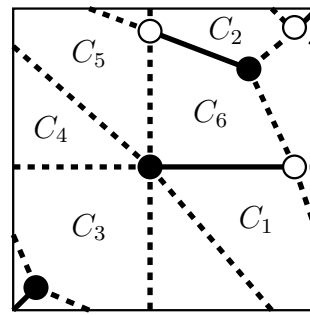

Figure 9.18: The dimer model and the perfect matching

Kazushi Ueda

Department of Mathematics, Graduate School of Science, Osaka University, Machikaneyama 1-1, Toyonaka, Osaka, 560-0043, Japan.

*e-mail address* :  kazushi@math.sci.osaka-u.ac.jp

Masahito Yamazaki

Department of Physics, Graduate School of Science, University of Tokyo, Hongo 7-3-1, Bunkyo-ku, Tokyo, 113-0033, Japan

*e-mail address* : yamazaki@hep-th.phys.s.u-tokyo.ac.jp